%% file: main.tex
\documentclass[letterpaper, 10 pt, journal, twoside]{ieeetran}

\input{header.tex}

\title{\LARGE \bf
Counter-examples in first-order optimization: a constructive approach}

\author{Baptiste Goujaud and Aymeric Dieuleveut and Adrien Taylor
\thanks{B. Goujaud is PhD candidate at
        CMAP, Ecole Polytechnique, Institut Polytechnique de Paris, Route de Saclay, 91120 Palaiseau, France.
        {\tt\small baptiste.goujaud@gmail.com}}
\thanks{A. Dieuleveut is faculty at CMAP,            Ecole Polytechnique, Institut Polytechnique de Paris,  Route de            Saclay, 91120 Palaiseau, France.
        {\tt\small aymeric.dieuleveut@polytechnique.edu}}
\thanks{A. Taylor is faculty at INRIA Paris,        2 Rue Simone IFF, 75012 Paris. {\tt\small adrien.taylor@inria.fr}}
}

\begin{document}

\maketitle
\thispagestyle{empty}
\pagestyle{empty}

\begin{abstract}
    While many approaches were developed for obtaining worst-case complexity bounds for first-order optimization methods in the last years, there remain theoretical gaps in cases where no such bound can be found. In such cases, it is often unclear whether no such bound exists (e.g., because the algorithm might fail to systematically converge) or simply if the current techniques do not allow finding them.
    In this work, we propose an approach to automate the search for cyclic trajectories generated by first-order methods. This provides a constructive approach to show that no appropriate complexity bound exists, thereby complementing approaches providing sufficient conditions for convergence. Using this tool, we provide ranges of parameters for which the famous Polyak heavy-ball, Nesterov accelerated gradient, inexact gradient descent, and three-operator splitting algorithms fail to systematically converge, and show that it nicely complements existing tools searching for Lyapunov functions.
\end{abstract}

\section{Introduction}
    \label{sec:intro}
    
    In the last years, first-order optimization methods (or algorithms) have attracted a lot of attention due to their practical success in many applications, including in machine learning (see, e.g.,~\cite{bottou2007tradeoffs}). Theoretical foundations for those methods played a crucial role in this success, e.g., by enabling the development of momentum-type methods (see, e.g.,~\cite{polyak_gradient_1963, nesterov1983method}). Formally, we consider the optimization problem
    \begin{equation}
        x_\star \triangleq \arg\min_{x\in \mathbb R^d} f(x) \tag{OPT}\label{eq:opt}
    \end{equation}
    for a  function $f$ belonging to a class of functions $\mathcal F$ (e.g., the set of convex functions, or the set of strongly convex and smooth functions, etc.). Classical first-order optimization methods for solving this problem include \emph{gradient descent}~(GD), \emph{Nesterov accelerated gradient method}~(NAG)~\cite{nesterov1983method}, and the \emph{heavy-ball method}~(HB)~\cite{polyak_gradient_1963}.
    These families of algorithms are parametrized: for example, GD is parametrized by a step-size $\gamma$ and HB is parametrized by both a step-size~$\gamma$ and a momentum parameter~$\beta$. 
    We generically denote by $\mathcal A$ any such method, for a specific choice of its parameters.
    For a given class of function $\mathcal{F}$ and an algorithm $\mathcal{A}$, we typically aim at answering the question
    \begin{center}
        \fbox{\parbox{0.48\textwidth}{
            \begin{center}
                Does $\mathcal{A}$ converge on every function of $\mathcal{F}$ \\ to their respective minimum?
            \end{center}
        }}
    \end{center}
    
    Common examples of function classes $\mathcal F$ {include} the set $\mathcal F_{\mu,L}$ of $\mu$-strongly convex and $L$-smooth functions, and the set $\mathcal Q_{\mu,L}$ of $\mu$-strongly convex and $L$-smooth \textit{quadratic} functions, for $\mu, L\geq 0$.
    
    This type of analysis, requiring results to hold on every function of a given class $\mathcal{F}$ is commonly referred to as \emph{worst-case analysis} and is the most popular paradigm for the analysis of optimization algorithms, see, e.g.,~\cite{nesterov1983method,dvurechensky2021first,bubeck2015convex,d2021acceleration,chambolle2016introduction}. 
    In this context, a very successful technique for proving worst-case convergence consists in looking for a decreasing sequence (called \textit{Lyapunov sequence}~\cite{lyapunov1992general,Kalman1960a,Kalman1960b}) of expressions $V_t$ of the iterates $x_t$, i.e.~such that
    \begin{equation}
        \forall f \in \mathcal{F},~ \forall t,~ \forall x_t, \quad V_{t+1}((x_s)_{s\leq t+1}) \leq V_t((x_s)_{s \leq t}), \label{eq:lyap}
    \end{equation}
    where some quantity of interest is upper-bounded by~ $V_T((x_s)_{s \leq T})$ as $T$ goes to infinity.
    For instance, when studying GD with step-size $1/L$ on the class $\mathcal{F}_{0, L}$ of $L$-smooth convex functions, we prove that $\forall f \in\mathcal{F}_{0, L}$, $ (t+1)(f(x_{t+1}) - f(x_\star))+ \tfrac{1}{2}\|x_{t+1}-x_\star\|^2 \leq t (f(x_{t})-f(x_\star)) + \tfrac{1}{2}\|x_{t}-x_\star\|^2$. 
    Therefore, $V_t((x_s)_{s \leq t}) = t (f(x_{t})-f(x_\star)) + \tfrac{1}{2}\|x_{t}-x_\star\|^2$ defines a decreasing sequence, and $f(x_t)-f(x_\star) \le V_t((x_s)_{s \leq t})/t \le V_0(x_0)/t $, proving convergence of this method on this class of functions.
    
    Due to the simplicity of the underlying proofs, the Lyapunov approach is particularly popular, e.g., for NAG~\cite{nesterov1983method,beck2009fast}, and HB~\cite{ghadimi2014global}. See~\cite{bansal2019potential,d2021acceleration} for surveys on this topic.
    
    \textbf{Necessary condition for worst-case convergence.}~While finding  a decreasing Lyapunov sequence guarantees convergence, not finding one does not guarantee anything: there may still exist a Lyapunov sequence, that the current analysis was not able to capture, or the method could converge without the existence of such Lyapunov sequence.
    Establishing that a method
    \emph{provably does not admit a worst-case convergence analysis} is therefore critical for avoiding spending an indefinite amount of time and effort searching for a non-existent convergence guarantee. 
    The existence of a \textit{cycle} for the algorithm on a particular function  means that it diverges on that function: in other words, the absence of cycle on all functions is a necessary condition for worst-case convergence.
    Moreover, a cycle can be observed  after only a \textit{finite} number of steps of the algorithm, while observing the divergence of a non periodic sequence is difficult or impossible. Overall,  this makes the search for cycles a computationally practical way of proving divergence.
    
    In order to discover cycles, we rely on computer-assisted worst-case analysis. \emph{Performance estimation problems} (PEP~\cite{drori2014performance, taylor2017smooth}) provide  a systematic approach to obtain convergence guarantees, including the search for appropriate Lyapunov arguments. Some packages (especially \texttt{Pesto}~\cite{taylor2017performance} and \texttt{Pepit}~\cite{goujaud2022pepit})  automate these tasks. We formulate cycle discovery as a minimization problem that can be cast in a PEP, and rely on the \texttt{Pepit} package to solve it.
    
    {\textbf{Examples:}} We demonstrate the applicability of our method on several examples.
    In particular, the case of HB illustrates the potential of our methodology. In fact, the search for the step-size $\gamma$ and momentum $\beta$ parameters leading to the fastest worst-case convergence over $\mathcal F_{\mu, L}$ is still an open problem, and the existence of parameters resulting in an accelerated rate remains a lingering question.  
    Indeed,~\cite{lessard2016analysis} exhibits a smooth and strongly convex function on which HB cycles, for   parameters $\gamma$ and $\beta$ optimizing the worst-case guarantee on $\mathcal Q_{\mu,L}$. 
    On the other hand,~\cite{ghadimi2014global} obtains a worst-case convergence on $\mathcal F_{\mu, L}$ for other parameters, but without acceleration. Recently,~\cite{upadhyaya2023automated} proposes a very general procedure to find Lyapunov sequences and extended the region of parameters $\gamma$ and $\beta$ HB provably converges on, leveraging PEPs. However, outside this region of the parameter space, the question of the convergence of the HB method remains open in the absence of a proof of divergence. For this example,  our approach demonstrates that a cycle exists for almost all parameters for which no Lyapunov is known.
    
    {\textbf{Summary of contributions:}}
    This paper proposes a systematic approach to prove that no worst-case certificate of convergence can be obtained for a given algorithm $\mathcal{A}$ on a class~$\mathcal{F}$. To do so, we establish the existence of a function in $\mathcal{F}$ over which $\mathcal{A}$ cycles.
    We illustrate our approach by applying it to three famous first-order  optimization algorithms, namely HB, NAG, inexact gradient descent with relatively bounded error. We further showcase the applicability of the approach to more general types of problems by studying the three-operator splitting method for monotone inclusions. For each method, we describe the set of parameters for which it is  known to converge and the ones where we establish the existence of a cycle. In the first three examples, our approach enables to fill the gap: we show the existence of cycles for all parametrizations not known to result in convergence.
    
    {\textbf{Organization:}}
    The rest of the paper is organized as follows. In \Cref{sec:definitions}, we introduce the concept of a stationary algorithm and formally define a cycle. In \Cref{sec:searching_for_cycles}, we present our methodology to discover cycles, relying on PEP. Finally, in \Cref{sec:examples}, we provide the numerical results.
    
\section{Definitions and notations}
    \label{sec:definitions}

    \noindent
    \begin{table}[t]
        \centering
        \begin{tabular}{ll}
        \toprule
        Notation & Corresponding object \\
        \midrule
            $\mathcal A$ & Generic algorithm \\
            $A$ & Update function of the algorithm $\mathcal A$ \\
            \eqref{eq:hb} & Heavy-ball \\
            \eqref{eq:nag} & Nesterov accelerated gradient \\
            \eqref{eq:igd} & Inexact GD\\
            \eqref{eq:tos} & Three operator splitting\\
            $\beta, \gamma $ & Algorithm parameters\\
            $(x_t)_t$ & Sequence of iterates generated by $\mathcal A$ \\
            $x_*$ & Optimal point \\
            $V$ & Lyapunov function \\
            $f$ & Objective function \\
            $\mathcal F$ & Generic class of functions \\
            $\mathcal F_{\mu, L}$ & Class of $L$-smooth and  $\mu$-strongly convex \\& functions \\
            $\mathcal Q_{\mu, L}$ & Class of  $L$-smooth and  $\mu$-strongly  convex \\ &quadratic functions \\
            $\ell $& Order of the algorithm $\mathcal A$ \\
            $ K $ & Length of the considered cycle \\
            $\mathcal O^{(f)} $ & Generic oracle applied on $f$ \\
            $u, F, G$  & Linearization variables (after SPD lifting)\\
            $d$ & Dimension\\
            $s_K$ & Score \\
            \bottomrule
        \end{tabular}
    \end{table}
    
    In this section, we consider a subclass of first-order methods, tailored for our analysis. It is chosen to ensure the periodicity of an algorithm that cycles once (see \Cref{prop:finite_horizon_cycle}). The class reduces to ``\emph{p-stationary canonical linear iterative optimization algorithms}'' (p-SCLI, see~\cite[Definition~1]{arjevani2016lower}) when the dependency to the previous iterates and gradients is linear which is a particular case of ``\emph{fixed-step first-order methods}'' (\emph{FSFOM}, see in~\cite[Definition~4]{taylor2017smooth}).
    Here, we consider \emph{stationary first-order methods (SFOM)}, whose iterates are defined as a fixed function of a given number of lastly observed iterates, as well as output of some oracles called on those iterates. 
    Examples of such oracles include gradients, approximate gradients, function evaluations, proximal step, exact line-search, Frank-Wolfe-type steps (see~\cite{taylor2017exact,goujaud2022pepit} for lists of oracles that can be handled using PEPs).
    The oracles we use depend on the setting under consideration.
    
    \begin{leftbot}
    \begin{Def}[Stationary first-order method (SFOM)]
    \label{def:SFOM}
        A method $\mathcal{A}$ is called order-$\ell$ \emph{stationary first-order  method} if there exists a deterministic first-order oracle $\mathcal{O}^{(f)}$ and a function $A$ such that the sequence generated on the function $f$ verifies $\forall t\geq\ell$,
        \begin{equation}
            x_{t} = A((x_{t-s}, \mathcal{O}^{(f)}(x_{t-s}))_{s\in\llbracket 1, \ell \rrbracket}). \tag{SFOM} \label{eq:sfom}
        \end{equation}
    \end{Def}
    \end{leftbot}
    
    For any given function of interest $f$ and any initialization $(x_t)_{t\in\llbracket 0, \ell -1 \rrbracket}$, an order-$\ell$ \eqref{eq:sfom} $\mathcal{A}$ iteratively generates a sequence $(x_t)_{t\in\mathbb{N}}$ that we denote $\mathcal{A}(f, (x_t)_{t\in\llbracket 0, \ell -1 \rrbracket})$.
    
    Definition~\ref{def:SFOM} above is very similar to the definition of a general first-order  method. However, the key assumption here is that the operation $A$ does not depend on the iteration counter~$t$: the algorithm is \textit{stationary}.
    While this assumption is restrictive, many first-order  methods are of the form \eqref{eq:sfom}, including (but not limited to): GD, HB~\cite{polyak_gradient_1963} and NAG~\cite{Nest03a} with constant step-sizes. On the other hand, any strategy involving decreasing step-size (e.g. for GD), or increasing momentum parameter (e.g. for NAG on $\mathcal F_{0,L}$ as in \cite{nesterov1983method}) are not in the scope of this definition.
    Note that the aforementioned examples use the first-order oracle $\mathcal{O}^{(f)}(x)\triangleq (\nabla f(x), f(x))$, although our methodology applies beyond this simple setting, as previously discussed. As an example, \ref{subsec:tos} considers an algorithm relying on the resolvent (or proximal operation).
    
    Stationarity is essential for being able to prove existence of a cyclical behavior in a finite number of steps. Next, we define a cyclic sequence.
    
    \begin{leftbot}
    \begin{Def}[Cycle]
        For any positive integer $K\geq2$, a sequence $(x_t)_{t \geq 0}$ is said to be $K$-cyclic if $\forall t \geq 0, x_t = x_{t+K}$.
        A sequence $x$ is said to be cyclic if there exists $K\geq2$ such that $x$ is $K$-cyclic.
    \end{Def}
    \end{leftbot}
    For any given order-$\ell$ \eqref{eq:sfom} $\mathcal{A}$, and any function class $\mathcal{F}$, we want to address the question
    \vspace{-0.5em}
    \begin{center}
        \fbox{\parbox{0.48\textwidth}{
            \begin{center}
                Does there exist a function $f \in \mathcal{F}$ and an initialization $(x_t)_{t\in\llbracket 0, \ell -1 \rrbracket}$ such that $\mathcal{A}(f, (x_t)_{t\in\llbracket 0, \ell -1 \rrbracket})$ is cyclic?
            \end{center}
        }}
    \end{center}
    
    \begin{Ex}
        \textit{
        In~\cite[Equation 4.11]{lessard2016analysis}, {the authors} answer positively to this question by providing a cycle of length~3, on the class $\mathcal{F}_{\mu, L}$ with $(\mu,L)=(1,25)$, and for $\mathcal{A}$ the heavy-ball method with step-size $\gamma =  (\frac{2}{\sqrt{L} + \sqrt{\mu}})^2$ and momentum parameter $\beta = (\frac{\sqrt{L} - \sqrt{\mu}}{\sqrt{L} + \sqrt{\mu}})^2$. Those parameters are natural candidates, that correspond to the limit of the step-size and momentum in Chebychev acceleration~\cite{flanders1950numerical,lanczos1952solution,young1953richardson}, and result in an acceleration for quadratic functions.
        }
    \end{Ex}
    
    In \Cref{sec:examples}, we extend this result to more parameters.

\section{Searching for cycles}
    \label{sec:searching_for_cycles}
    
    In this section, we show how to find  cyclic trajectories.
    
    \subsection{Motivation}
        Finding diverging trajectories for an algorithm $\mathcal A$ might be challenging. We thus focus on cycles, as they allow to focus on a finite sequences of iterates only. Indeed, for an SFOM,  once we observe the cycle to be repeated once, we can easily extrapolate: this same cycle is repeated again and again. This statement is formalized in the following proposition.
        
        \begin{leftbot}
        \begin{Prop}
            \label{prop:finite_horizon_cycle}
            Let $\mathcal{A}$ be a order-$\ell$ \eqref{eq:sfom}, and $(x_t)_{t\in\mathbb{N}}$ be any sequence generated by $\mathcal{A}$.
            Then the sequence $(x_t)_{t\in\mathbb{N}}$ is cyclic if and only if there exists $K\geq2$ such that $\forall t \in \llbracket 0, \ell-1 \rrbracket, x_t = x_{t+K}$.
        \end{Prop}
        \end{leftbot}

        \begin{proof}{}
            Let $\mathcal{A}$ be a order-$\ell$ \eqref{eq:sfom}, and $(x_t)_{t\in\mathbb{N}}$ be any sequence generated by $\mathcal{A}$. The method $\mathcal A$ is cyclic if and only if there exists $K\geq 2$, such that  the translated sequence $(\Tilde{x}_t)_{t\in\mathbb{N}}:=  (x_{t+K})_{t\in\mathbb{N}}$, is identical to $({x}_t)_{t\in\mathbb{N}}$. \Cref{prop:finite_horizon_cycle} states that those two sequences are identical if and only if their $\ell$ first terms are. It is clear that if the sequences $x$ and $\Tilde{x}$ are identical, their $\ell$ first terms also are. Reciprocally, let's assume that their $\ell$ first terms are identical and let's introduce the function $f$ and associated oracles $\mathcal{O}^{(f)}$ defined such that $x = \mathcal{A}(f, (x_t)_{t\in\llbracket 0, \ell -1 \rrbracket})$.
            Then $\forall t \geq 0, x_{t} = A((x_{t-s}, \mathcal{O}^{(f)}(x_{t-s}))_{s\in\llbracket 1, \ell \rrbracket})$.
            In particular  $\forall t \geq 0$
            \begin{align*}
                x_{t+K} & = A((x_{t+K-s}, \mathcal{O}^{(f)}(x_{t+K-s}))_{s\in\llbracket 1, \ell \rrbracket}),
            \end{align*}
            thereby reaching
            \begin{align*}
            \Tilde{x}_{t} & = A((\Tilde{x}_{t-s}, \mathcal{O}^{(f)}(\Tilde{x}_{t-s}))_{s\in\llbracket 1, \ell \rrbracket}).
            \end{align*}
            
            Consequently, $\tilde x = \mathcal{A}(f, (\Tilde{x}_t)_{t\in\llbracket 0, \ell -1 \rrbracket}) $, and since the $\ell$ first terms of $x$ and $\Tilde{x}$ are identical, \\
            $\Tilde{x} = \mathcal{A}(f, (\Tilde{x}_t)_{t\in\llbracket 0, \ell -1 \rrbracket}) = \mathcal{A}(f, (x_t)_{t\in\llbracket 0, \ell -1 \rrbracket}) = x$.
        \end{proof}

    \subsection{Approach}
        We now present the approach used to search for cycles, based on performance estimation problems (PEPs)~\cite{drori2014performance,taylor2017smooth}. We consider an algorithm $\mathcal{A}$, a function $f$ and initial points $(x_t)_{t\in\llbracket 0, \ell -1 \rrbracket}$, and run $\mathcal{A}$ on $f$ starting on $(x_t)_{t\in\llbracket 0, \ell -1 \rrbracket}$. This generates the sequence $x = \mathcal{A}(f, (x_t)_{t\in\llbracket 0, \ell -1 \rrbracket})$. For any positive integer~$K$, we then define the non-negative score 
        \begin{equation*}
            s_K(\mathcal{A}, f, (x_t)_{t\in\llbracket 0, \ell -1 \rrbracket}) = \sum_{t=0}^{\ell -1} \|x_t - x_{t+K}\|^2.
        \end{equation*}
        From \Cref{prop:finite_horizon_cycle}, this score is identically zero if and only if $\mathcal{A}$ cycles on $f$ when starting from $(x_t)_{t\in\llbracket 0, \ell -1 \rrbracket}$.
        This suggests that one can search for cycles of length $K$ by minimizing the score $s_K(\mathcal{A}, f, (x_t)_{t\in\llbracket 0, \ell -1 \rrbracket})$ w.r.t.~the function  $f$ and the initialization $(x_t)_{t\in\llbracket 0, \ell -1 \rrbracket}$.
        
        Observe that fixed points of $\mathcal{A}$, that correspond to cycles of length 1, also cancel this score. Our goal is to search for cycles of length at least $K\geq2$, that entail that the algorithm diverges for a particular function and initialization. As any convergent algorithm must admit the optimizer of $f$ as fixed point, we have to exclude fixed points. To do so, we add the constraint that the two first iterates are far from each other. In most cases of interest, making this constraint can be done without loss of generality due to the homogeneity of the underlying problems.
        We arrive to the following formulation:

        \begin{equation}
            \left|
            \begin{array}{cc}
                \underset{d \geq 1, f \in\mathcal{F}, x\in\left(\mathbb{R}^d\right)^{\mathbb{N}}}{\text{minimize }} & \sum_{t=0}^{\ell -1} \|x_t - x_{t+K}\|^2 \\
                \text{subject to } & 
                \left\{
                \begin{array}{c}
                    x = \mathcal{A}(f, (x_t)_{t\in\llbracket 0, \ell -1\rrbracket}) \\
                    \|x_1 - x_0\|^2 \geq 1.
                \end{array}
                \right.
            \end{array}
            \right.
            \tag{$\mathcal{P}$} \label{eq:problem}
        \end{equation}
        
        As we see in the next sections, this problem can be used to answer the question of interest by testing the nullity of the solution of \eqref{eq:problem}.

        As is, \eqref{eq:problem} looks intractable due to the minimization over the infinite-dimensional space $\mathcal{F}$ and its non-convexity. This can be handled using the techniques proposed in~\cite{taylor2017smooth,  taylor2017exact}, developed for PEP. It consists in reformulating \eqref{eq:problem} into a semi-definite program (SDP) using interpolation / extension properties for the class $\mathcal{F}$, together with SDP lifting.

        Indeed, \eqref{eq:problem} does not fully depend on $f$, but only on $\mathcal{O}^{(f)}(x_t)$ where $t \in \llbracket 0 ;K+\ell-2\rrbracket$.
        By introducing the variables $\mathcal{O}_{ t} \triangleq \mathcal{O}^{(f)}(x_t)$, we can replace the constraint $x = \mathcal{A}(f, (x_t)_{t\in\llbracket 0, \ell -1\rrbracket})$ of \eqref{eq:problem} by
        \begin{equation*}
            \left\{
            \begin{array}{ccll}
                x_{\ell} & = & A(( & \hspace{-.3cm} x_{\ell -s}, \mathcal{O}_{ \ell -s})_{s\in\llbracket 1, \ell \rrbracket}), \\
                & \vdots && \\
                x_{K + \ell - 1} & = & A(( & \hspace{-.3cm} x_{K+\ell - 1-s}, \mathcal{O}_{ K+\ell - 1-s})_{s\in\llbracket 1, \ell \rrbracket}),
            \end{array}
            \right.
        \end{equation*}
        
        \noindent and minimize over the finite dimensional variables $(\mathcal{O}_{t})_{0 \leq t \leq K+\ell-2}$ instead of $f$, under the constraint that there exists a function $f\in\mathcal{F}$ that interpolates those values, i.e.~that verifies $\mathcal{O}^{(f)}(x_t) = \mathcal{O}_{t}$ for all $t \in \llbracket 0 ;K+\ell-2\rrbracket$. For some classes $\mathcal F$, those  interpolation property are equivalent to tractable inequalities, as in the following example.

        \begin{Ex}[$L$-smooth convex functions]\label{ex:ic_smooth}
        \textit{
        If the oracles are only the gradients and the function values of the objective function $f$, denoting $f_i \triangleq f(x_i)$ and $g_i \triangleq \nabla f(x_i)$  (i.e.~$\mathcal{O}_{i}\triangleq (g_i, f_i)$), the interpolation conditions of $\mathcal{F}_{0, L}$ are provided in~\cite{taylor2017smooth} as
        \begin{equation}
            \forall i, j,~ f_i \geq f_j + \left< g_j, x_i - x_j \right> + \tfrac{1}{2L}\|g_i - g_j\|^2. \tag{IC} \label{eq:cni}
        \end{equation}
        }
        \end{Ex}
        This function class is considered in three of the four examples under consideration in the next section (HB, NAG, inexact GD).
        However, the methodology described in this paper applies to many other classes beyond $\mathcal F_{\mu,L}$, see for instance \cite[\href{https://pepit.readthedocs.io/en/0.2.1/api/functions_and_operators.html}{Function classes}]{goujaud2022pepit}: an example of such a class is used in the fourth example of the following section. Each class considered must be described by its interpolation conditions, similar to \eqref{eq:cni}. Other examples of known interpolations conditions are provided in \cite[Th. 3.4-3.6]{taylor2017exact}, \cite[Cor.1\&2]{dragomir2021optimal}, \cite[Th.1]{guille2022gradient} \cite[Th. 2.6]{goujaud2022optimal}.
        The key ingredient for a class to enter in the scope of this paper, is that its interpolation conditions are expressed as a degree~2 polynomial in $x_t$ and $\mathcal{O}_t$, and that a given variable is not involved both in a monomial of degree 1 and one of degree~2, as in \eqref{eq:cni}.
        
        Then, the SDP lifting part consists in introducing a Gram matrix~$G$ \cite[Theorem~5]{taylor2017smooth} of vectors among $x_t$ and $\mathcal{O}_t$ that are involved in degree 2 monomials, so that those quadratic expressions of $x_t$ and $\mathcal{O}_t$ are then expressed linearly in term of $G\succeq 0$. Thereby, the problem can be cast a standard SDP.

        In the case where the oracle is $\mathcal{O}_t = (g_t, f_t)$, 
        and the class of interest $\mathcal{F}$ is the class of $L$-smooth convex functions $\mathcal{F}_{0, L}$, the objective and all the constraints of \eqref{eq:problem} are written linearly in terms of $(f_t)_t$ and quadratically in terms of $(x_t, g_t)_t$. Therefore, we define $G$ as the Gram matrix of $(x_t, g_t)_t$ and $F$ as a vector storing all the values $f_t$ leading to an SDP reformulation of the problem. See, e.g.,~\cite[Section 2]{goujaud2022pepit} for a detailed derivation on a simple example.

        Setting $u \triangleq (G, F)$, \eqref{eq:problem} is generally rewritten, under above mentioned key ingredients, as

        \begin{equation}
            \left|
            \begin{array}{cc}
                \underset{u}{\text{minimize }} & \left< u, v_{obj} \right> \\
                \text{subject to } & 
                \left\{
                \begin{array}{c}
                    \left< u, v_1 \right> \geq 0 \\
                    \ldots \\
                    \left< u, v_n \right> \geq 0 \\
                    \left< u, v_{\text{aff}} \right> \geq 1 \\
                    u \in \mathcal{C}.
                \end{array}
                \right.
            \end{array}
            \right.
            \tag{SDP-$\mathcal{P}$} \label{eq:sdp_problem}
        \end{equation}

        The objective is linear, as well as the  first $n$ constraints.
        The affine constraint  $\left< u, v_{\text{aff}} \right> \geq 1 $ enables to discard the trivial solution $u=0$ and corresponds in \eqref{eq:problem} to the constraint $\|x_1 - x_0\|^2 \geq 1$.
        Finally, the constraint $u\in\mathcal{C}$ corresponds to the constraint $G \succeq 0$. $\mathcal{C}$ is then a closed convex semi-cone.

        By definition, if there exists a feasible vector $u$ such that the objective of \eqref{eq:sdp_problem} is zero, then it describes a cycle.
        Moreover \eqref{eq:sdp_problem} is convex and efficiently solvable (due to the existence of a Slater point~\cite[Theorem 6]{taylor2017smooth}).
        
        In the next sections, we numerically apply this methodology through the \texttt{Pepit} python package~\cite{goujaud2022pepit} which takes care about the tractable reformulations of \eqref{eq:problem} into~\eqref{eq:sdp_problem}. \eqref{eq:sdp_problem} is then solved using through a standard solver~\cite{mosek} to determine the infimum value of $\left< u, v_{\text{obj}} \right>$ over the feasible set of \eqref{eq:sdp_problem}. The next theorem allows to conclude about the existence of cycles.

        \begin{leftbot}
        \begin{Th}
            \label{thm:min=inf}
            Assuming the infimum value of \eqref{eq:problem} to be 0, then there exists a cycle.
        \end{Th}
        \end{leftbot}
        
        \begin{proof}{}
            By equivalence between~\eqref{eq:problem} and~\eqref{eq:sdp_problem}, we assume that there exists a sequence of feasible vectors $u_i$ with $\langle u_i, v_{\text{obj}} \rangle \rightarrow 0$.
        
            The constraint $\left< u_i, v_{\text{aff}} \right> \geq 1$, guarantees that none of the~$u_i$ is equal to 0. 
            Considering any norm (all equivalent to each other in finite dimension) and projecting $u_i$ on the associated sphere defines $s_i \triangleq \tfrac{u_i}{\|u_i\|}$. The other $n$ \textit{linear} constraints still hold after the scaling and $s_i \in \mathcal{C}$ since $\mathcal{C}$ is a semi-cone.
            
            Moreover, the norms of all $u_i$ are lower bounded by the distance from 0 to the affine hyperplan $\left\{ w | \left< w, v_{\text{aff}} \right> = 1 \right\}$. Hence $|\left< s_i, v_{\text{obj}} \right>| = \tfrac{|\left< u_i, v_{\text{obj}} \right>|}{\|u_i\|} \rightarrow 0$.
        
            Finally by compacity of the sphere, there exists a subsequent limit $s$ of the sequence $(s_i)_i$ and by continuity of the linear operator $\left< \cdot , v_{\text{obj}} \right>$, $\left< s , v_{\text{obj}} \right>=0$.
        
            We conclude that $s$ is a vector the objective reaches 0 on, that verifies all the constraints of \eqref{eq:sdp_problem} but the affine one.
        
            Note the affine constraint only aimed at discarding the trivial solution 0 in a linear way (for solver purpose), and that $s$ is not 0. Then $s$ describes a cycle.
        \end{proof}

\section{Application to four different \texorpdfstring{\eqref{eq:sfom}s}{SFOMs}}
    \label{sec:examples}
    In this section we illustrate this methodology on four examples: heavy-ball (HB), Nesterov accelerated gradient (NAG), gradient descent (GD) with inexact gradients, and three-operator splitting (TOS). For each, we apply the methodology proposed in \Cref{sec:searching_for_cycles}. The code is available in the public github repository \url{https://github.com/bgoujaud/cycles}. Since ingredients are the same as those of classical PEPs, we also use the python package \texttt{Pepit}~\cite{goujaud2022pepit}. We perform a grid search over the spaces of \textit{parameters of interest}~$\Omega$, described in the respective subsections.
    We compare the parameter region where Lyapunov functions can be obtained with the region in which we establish that the method cannot have a guaranteed worst-case convergence (due to the existence of cycles). More precisely, in \Cref{fig:hb,fig:nag,fig:igd,fig:tos} below,  green regions corresponds to parameter choices for which the methods converge (existence of a Lyapunov function, found using the code provided with~\cite{taylor2018lyapunov}). Conversely, in the red regions, our methodology establishes that the method cycles on at least one function of $\mathcal{F}$.
    In short, the algorithms converge in the green regions and do not converge in the worst-case in the red ones.
    
    Note that some parameters for which the algorithm $\mathcal A$ admits a worst-case convergence guarantee could theoretically exist outside the green region: indeed, in~\cite{taylor2018lyapunov}, the authors do not guarantee that they necessarily find convergence.
    Similarly, parameters for which $\mathcal A$ does \textit{not} admit a worst-case convergence guarantee could theoretically exist outside of the red region: indeed
    \eqref{eq:problem} is defined for a fixed cycle length $K$, and we therefore run it several times with different values of~$K$. Longer cycles are therefore not detected. Moreover, the non-existence of cycles does not necessarily imply that the algorithm always converges.
    
    Interestingly, in practice,  we observe on~\Cref{fig:nag,fig:igd} that the set $\Omega$ of parameters of interest is completely filled by the union of those 2 regions and that it is almost the case on~\Cref{fig:hb} (it may have been if we had searched for cycles of all lengths). As a consequence, we fully characterize for which tuning the algorithms admit a guaranteed worst-case convergence. On the contrary, there remains a significant gap between the red and the green regions in our last example, see~\Cref{fig:tos}.
    
    \subsection{Heavy-ball}

        The HB algorithm, as introduced by~\cite{polyak_gradient_1963}, corresponds to the following update, for a step-size parameter $\gamma$ and a momentum parameter $\beta$:
        \begin{equation}
            x_{t+1} = x_t + \beta (x_t - x_{t-1}) - \gamma \nabla f(x_t). \tag{HB} \label{eq:hb}
        \end{equation}
        Therefore \eqref{eq:hb} is an order-2 \eqref{eq:sfom}.
        
        \begin{figure}[t]
            \centering
            \includegraphics[width=\linewidth]{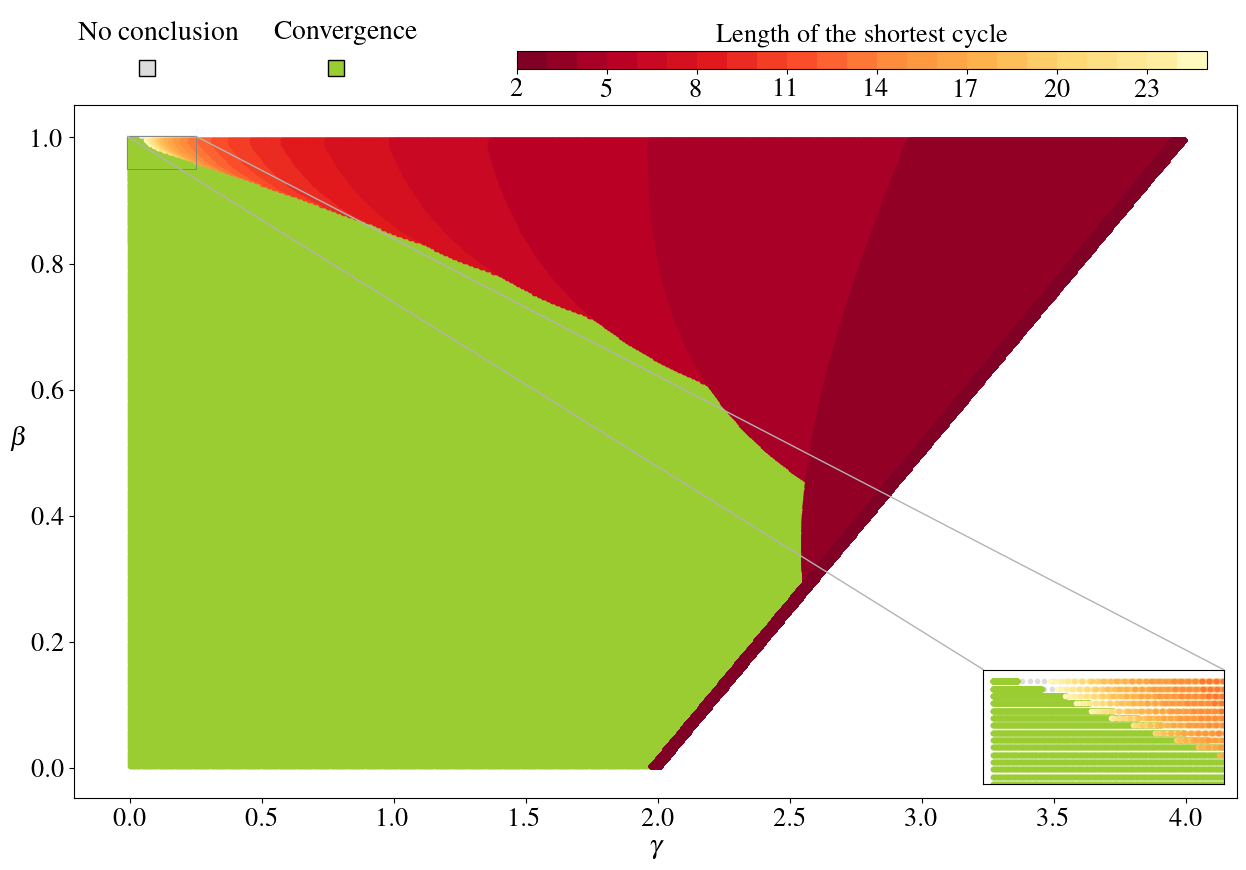}
            \vspace{-1em}\caption{Heavy-ball~\eqref{eq:hb}. Green area: set of parameters $(\gamma, \beta)\in \Omega_{\mathrm{HB}}$ for which a Lyapunov function exists; Red area: set of parameters $(\gamma, \beta) \in \Omega_{\mathrm{HB}}$ for which \eqref{eq:hb} cycles on at least one function in $\mathcal F_{0,L}$.}
            \label{fig:hb}
        \end{figure}
        
        \textbf{Set $\Omega_{\mathrm{HB}}$ of parameters of interest:} HB converges on the set  $\mathcal Q_{0,L}$ if and only if the parameters $\gamma, \beta$ verify $0 \leq \gamma \leq {2(1+\beta)}/{L} \leq {4}/{L}$~\cite{polyak_gradient_1963}. Note this condition enforces $-1 \leq \beta \leq 1$. Moreover, we restrict to $\beta \geq 0$ as $\beta<0$ is not an interesting setting (slowing down convergence with respect to GD).  Therefore, we limit our analysis to this set of parameters. 
        
        \textbf{Interpretation.} The red area in \Cref{fig:hb} shows parameters where cycles of length $K \in \llbracket2; 25\rrbracket$ are found by our methodology. 
        The red color intensity indicates the length of the shortest cycle.
        
        A striking observation is that the space $\Omega_{\mathrm{HB}}$ is almost filled by the union of the red area and green one (where Lyapunov functions exist). Thereby, for almost all values of the parameters, we have a definitive answer on the existence of a certificate of convergence in the worst-case. That being said, there exists a small unfilled region in the top left corner (see the zoom on \Cref{fig:hb}) In this region, we do not know how HB behaves, and whether it accelerates. However, adding longer cycle length may enable to obtain cycles in that area. Indeed we considered only cycles of length $K \le 25$, for computational reasons.
        Recently,~\cite{goujaud2023provable} computed the analytical region of cycles of~\eqref{eq:hb} showing non-acceleration of the latter.

    \subsection{Nesterov accelerated gradient}

        \begin{figure}[t]
            \centering
            \includegraphics[width=\linewidth]{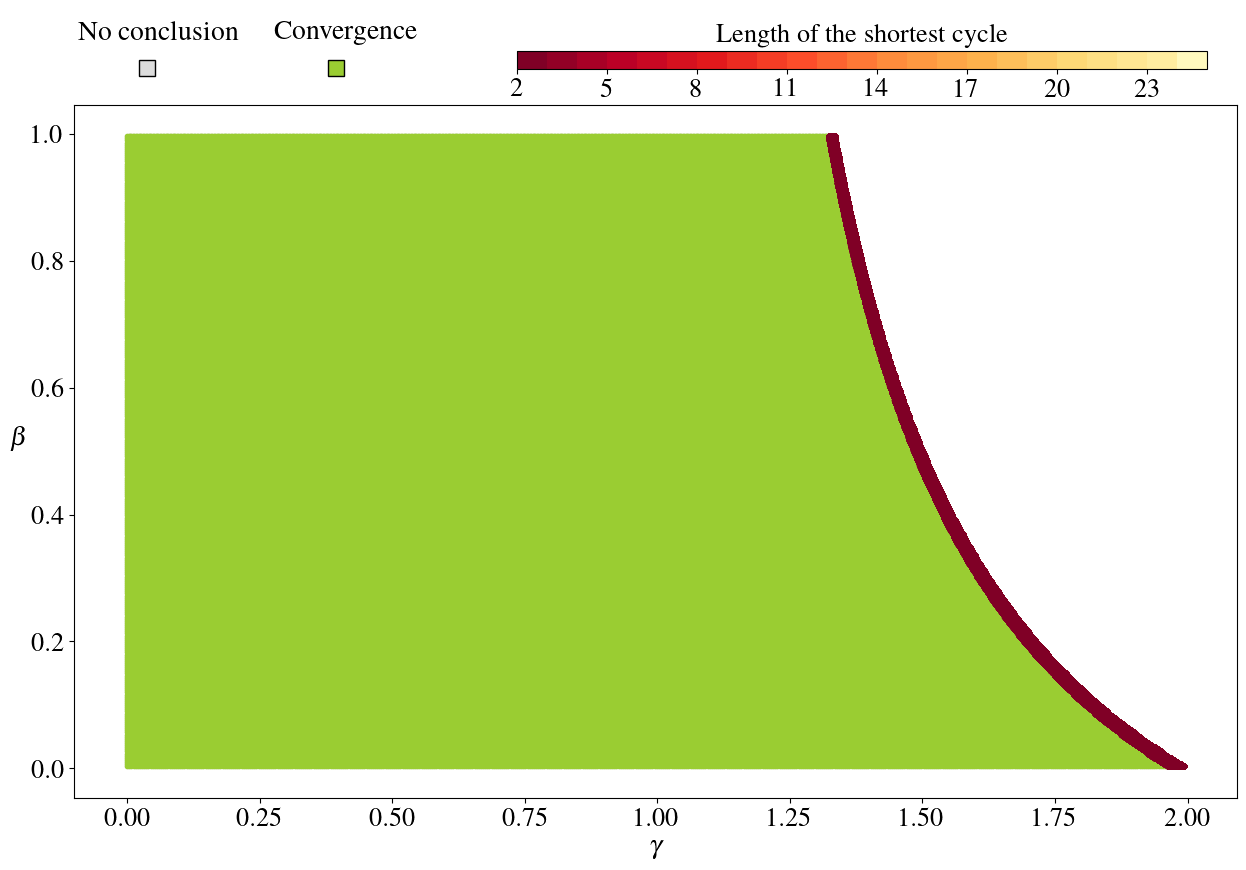}
            \vspace{-1em}\caption{Nesterov Accelerated gradient~\eqref{eq:nag}. Green area: set of parameters $(\gamma,\beta)\in \Omega_{\mathrm{NAG}}$ for which a Lyapunov function exists; Red area: set of $( \gamma, \beta) \in \Omega_{\mathrm{NAG}}$ for which \eqref{eq:nag} cycles on at least one function in $\mathcal F_{0,L}$.} \label{fig:nag}
        \end{figure}

        NAG (also known as the \textit{fast gradient method}) was introduced by~\cite{Nest03a} and corresponds to the following update, for a step-size parameter~$\gamma$ and a momentum parameter $\beta$:
        \begin{equation}
           \left \lbrace  \begin{array}{ccl}
                y_t &= & x_t + \beta (x_t - x_{t-1}), \\
                x_{t+1} & =& y_t - \gamma \nabla f(y_t).
            \end{array}\right.
            \tag{NAG} \label{eq:nag}
        \end{equation}
        \eqref{eq:nag} is also written as follows
        $y_{t+1} = (1 + \beta) (y_t - \gamma \nabla f(y_t)) - \beta (y_{t-1} - \gamma \nabla f(y_{t-1}))$,
        and is therefore also an order-2 \eqref{eq:sfom}.
        
        \textbf{Set $\Omega_{\mathrm{NAG}}$ of parameters of interest:} As for HB, we consider the set of $\beta$ and $\gamma$  for which~\eqref{eq:nag} converges on $\mathcal Q_{0,L}$. This corresponds to considering all $\beta, \gamma$ verifying $0\leq\beta\leq 1$ and $0\leq{\gamma}\leq \tfrac{2}{L}\tfrac{1 + \beta}{1+2\beta}$.
        
        \textbf{Interpretation:} \eqref{eq:nag} is known to converge, with an accelerated rate, on $\mathcal F_{\mu,L}$,  for the  tuning $(\gamma,\beta)= (\frac{1}{L}, \tfrac{\sqrt{L} - \sqrt{\mu}}{\sqrt{L} + \sqrt{\mu}})$, that optimizes the convergence rate on $\mathcal Q_{\mu,L}$. For this reason,~\eqref{eq:nag} is considered to be more ``robust'' than HB.
        
        \Cref{fig:nag} shows that \eqref{eq:nag} admits a Lyapunov function for almost any parameters in $\Omega_{\mathrm{NAG}}$.
        Moreover, our methodology does not detect any set of parameters at which a cycle  of length $K \in \llbracket 2; 25\rrbracket$  exists, apart on the boundary $\{(\gamma, \beta), \gamma= \tfrac{2}{L}\tfrac{1 + \beta}{1+2\beta}\}$. On the boundary, cycles of length 2 are observed, whose existences are theoretically  verified  on one-dimensional quadratic functions. This illustrates the robustness of our methodology when very few cycles exist.

    \subsection{Inexact gradient method}
        
        Next, we consider the inexact gradient method, parameterized by $\gamma$ and $\varepsilon$ and the update
        \begin{equation}
            \begin{array}{rl}
                \text{Get } \mathcal{O}(x_t)\ = &  d_t  \ \text{ such that } \|d_t - \nabla f (x_t)\| \leq \varepsilon \|\nabla f (x_t)\|, \\
                x_{t+1} \  = & x_t - \gamma d_t.
            \end{array}
            \tag{IGD} \label{eq:igd}
        \end{equation}
        \eqref{eq:igd} is thus an~\eqref{eq:sfom} of order 1.

        \begin{figure}[t]
            \centering
            \includegraphics[width=\linewidth]{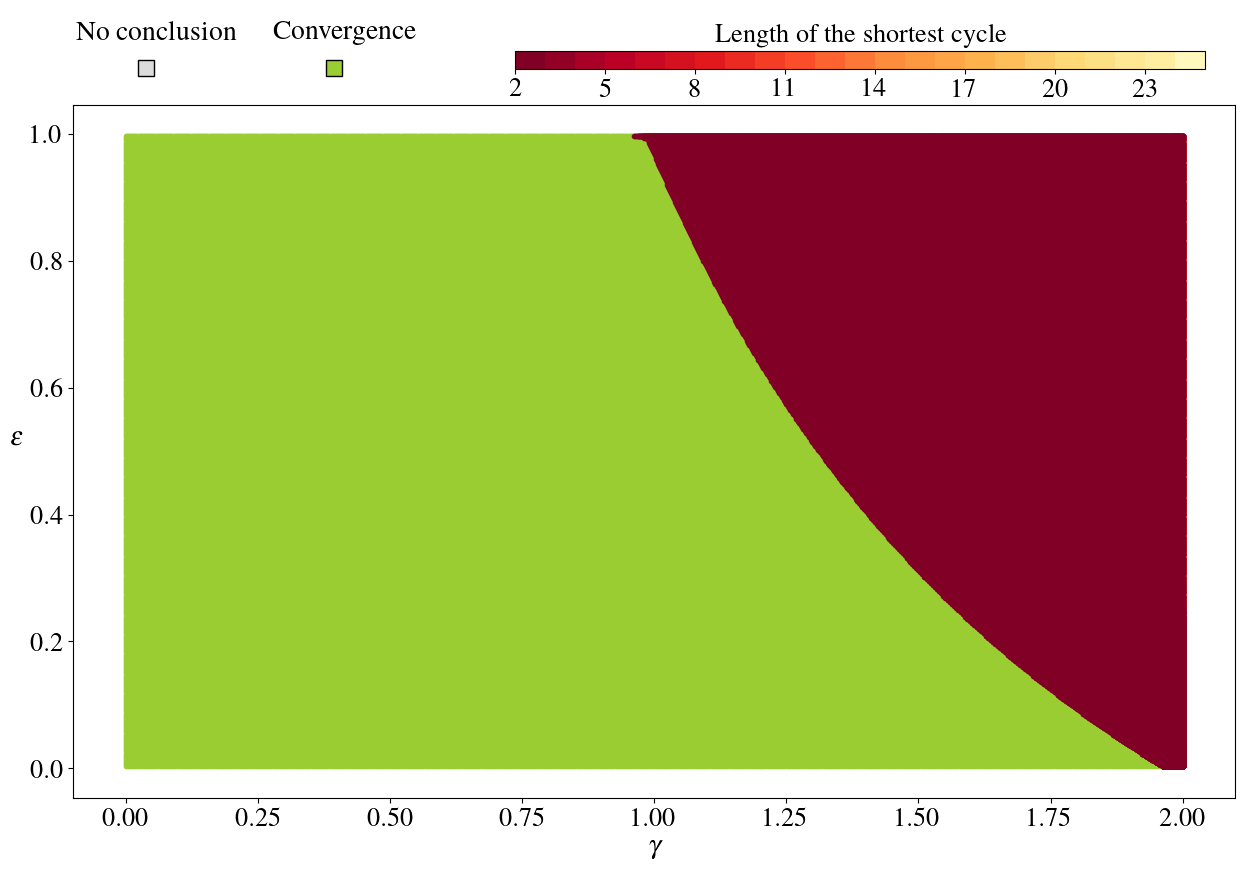}
            \vspace{-1em}\caption{Inexact GD~\eqref{eq:igd}. Green area: set of parameters $(\gamma,\varepsilon)\in \Omega_{\mathrm{IGD}}$ for which a Lyapunov function exists; Red area: set of $(\gamma,\varepsilon)\in \Omega_{\mathrm{IGD}}$ for which \eqref{eq:igd} cycles on at least one function in $\mathcal F_{0,L}$.}
            \label{fig:igd}
        \end{figure}
        
        \textbf{Set $\Omega_{\mathrm{IGD}}$ of parameters of interest:} Since the exact gradient method converges only for $\gamma < \frac{2}{L}$, we only consider such steps-sizes. Moreover, $\varepsilon \geq 1$ allows $d_t = 0$ and thereby does not make much sense. This motivates considering the set $\Omega_{\mathrm{IGD}}=\{(\gamma, \varepsilon) \in [0; \frac{2}{L}]\times [0,1]\}$.
        
        \textbf{Interpretation:}  We search for cycles of length $K \in \llbracket 2; 25\rrbracket$ and use color intensity to  show the minimal cycle length. \eqref{eq:igd} is known to converge for any $\gamma \leq \frac{2}{L(1 + \varepsilon)}$~(see~\cite{de2020worst,gannot2022frequency}). \Cref{fig:igd} shows that the complementary of this region of convergence is completely filled by parameters allowing cycles, showing that no other parameters values than the known ones allow obtaining worst-case convergence of \eqref{eq:igd}.
        
    \subsection{Three-operator splitting}\label{subsec:tos}

        \begin{figure}[b]
            \centering
            \includegraphics[width=\linewidth]{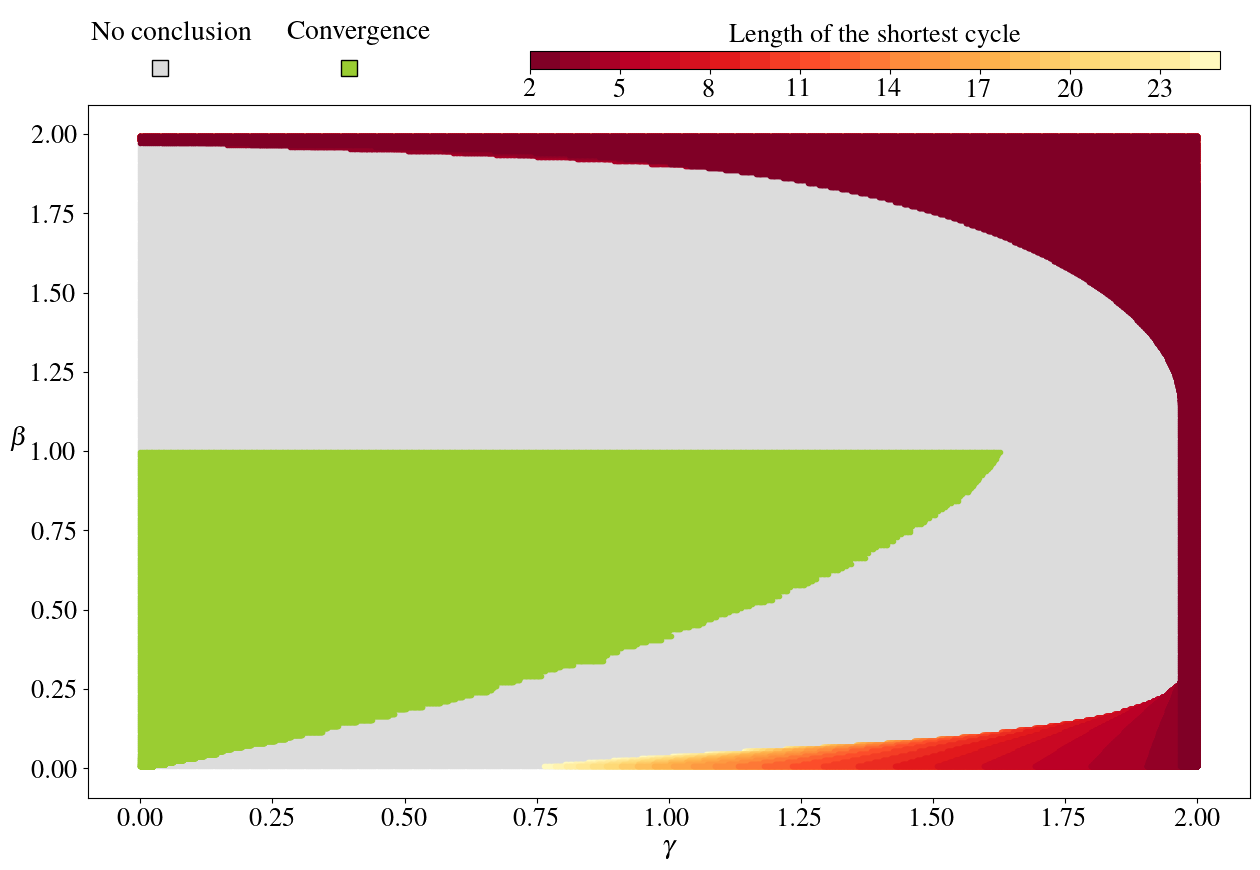}
            \vspace{-1em}\caption{Three-operator splitting~\eqref{eq:tos}. Green area: set of parameters $(\gamma,\beta)\in \Omega_{\mathrm{TOS}}$ for which a Lyapunov function exists; Red area: set of $(\gamma,\beta)\in \Omega_{\mathrm{TOS}}$ for which \eqref{eq:tos} cycles on at least a triplet of operators.}
            \label{fig:tos}
        \end{figure}

        The three-operator splitting (TOS) method, introduced by \cite{davis2017three}, aims at solving the \emph{inclusion problem} $0 \in Ax + Bx + \partial f(x)$, where $A$ is a monotone operator, $B$ is a co-coercive operator and $\partial f$ denotes the differential of the smooth (strongly) convex function $f$.
        It corresponds to the following update, for a step-size parameter $\gamma$, a smoothing parameter $\alpha$ and an update parameter $\beta$:
        \begin{equation}
            \left \lbrace
            \begin{array}{ccl}
                x_{t+1} & = & J_{\alpha B} (w_t), \\
                y_{t+1} & = & J_{\alpha A} \left(2 x_{t+1} - w_t - \frac{\gamma}{\beta} \nabla f(x_{t+1})\right), \\
                w_{t+1} & = & w_t - \beta (x_{t+1} - y_{t+1}),
            \end{array}\right.
            \tag{TOS} \label{eq:tos}
        \end{equation}
        where $J_{O}$ denotes the resolvent of the operator $O$, i.e. $J_O = (I + O)^{-1}$.
        Note that \eqref{eq:tos} is therefore an order-1 \eqref{eq:sfom}.
        
        \textbf{Set $\Omega_{\mathrm{TOS}}$ of parameters of interest:} When considering $A$, $B$ and $\nabla f$ to be linear symmetric and co-diagonalisable operators, the set of convergence of \eqref{eq:tos} is $\Omega_{\mathrm{TOS}} = \left\{ (\gamma, \beta) \in \left[0, \frac{2}{L}\right] \times [0, 2] \right\}$, which we therefore consider.
        
        \textbf{Interpretation:} We search for cycles of length $K\in\llbracket 2, 25\rrbracket$ and use color intensity to show the minimal cycle length. Interestingly, there is a gap between the green region and the red ones. Unlike for \eqref{eq:hb}, it seems that increasing the length of the cycle does not help covering this gap and shows that some algorithms might have no Lyapunov function while not cycling. Understanding the behavior of \eqref{eq:tos} in the grey region is therefore still an open question.

\section{Conclusions}
    \label{sec:conclusions}

    \textbf{Summary.}
    This work proposes a systematic approach for finding counter-examples to convergence of first-order methods, bringing a complementary tool to the existing systematic techniques for finding convergence guarantees (that include certifications through the existence of Lyapunov functions). Our approach is based on now classical tools and techniques used in the field of first-order optimization and a few existing packages~\cite{goujaud2022pepit,taylor2017performance} allows for straightforward implementations of our methodology.

    \textbf{Discussion and Future works.}
    While our analysis complements the Lyapunov one, existence of a Lyapunov function or existence of a cycle are not the only 2 options. In this section, we discuss the eventuality that an algorithm diverges on at least one function in a class, without resulting in cycles.

    Indeed, the sequence of iterates produced by an algorithm on a given function may diverge by (i)~tending to infinity, (ii)~simply growing unbounded, or even (iii)~showing a chaotic behavior, while staying in a compact set.
    Understanding if those three divergence cases might occur on functions $f\in \mathcal F$ while \textit{no} function  $f\in \mathcal F$ results in a cycle, is thus an interesting open question:
    \vspace{-1em}
    \begin{center}
        \fbox{\parbox{0.48\textwidth}{
            \begin{center}
                If the class $\mathcal{F}$ is connected, the update $A$ of an algorithm $\mathcal{A}$ and the oracle $\mathcal{O}$ are continuous, if $\mathcal{A}$ converges on one function of $\mathcal{F}$ and diverges on another one, is there necessarily one function of $\mathcal{F}$ which $\mathcal{A}$ cycles on?
            \end{center}
        }}
    \end{center}

    An interesting example is $\mathcal{A}$ being GD with step-size 1 on the ($1+\rho$)-smooth (and $2\rho$-Hessian-Lipschitz) function $f_\rho$, such that $f_\rho(x)$ is equal to:
    \begin{equation*}
       \left\{
        \begin{array}{ll}
            \tfrac{\rho}{3} |x|^3 + \tfrac{1-\rho}{2} x^2 + \frac{(\rho-1)^3}{6\rho^2} & \text{if } |x| \leq 1, \\
            -\tfrac{\rho}{3} |x|^3 + \tfrac{1+3\rho}{2} x^2 - 2 \rho |x|+ \frac{(\rho-1)^3 + 4\rho^3}{6\rho^2} & \text{if } 1 \leq |x| \leq \tfrac{3}{2}, \\
            \tfrac{1}{2}x^2 + \tfrac{\rho}{4} |x| + \frac{4(\rho-1)^3-11\rho^3}{24\rho^2} & \text{if } \tfrac{3}{2} \leq |x|,
        \end{array}
        \right.
    \end{equation*}
    for $\rho \in [0, 4]$.
    From any point in the interval $(1, \infty)$, the next iterate is in $[-1, 0]$.
    Similarly, starting in the interval $(-\infty, -1)$, the second iterate is in $[0, 1]$.
    Those 2 intervals are stable and, by symmetry of the function, the dynamics in those 2 intervals are themselves symmetric.
    Note that for any $x\in[0, 1]$, $f'(x) = \rho x^2 + (1-\rho)x$, leading to the dynamic $x_{t+1} = x_t- 1\times \nabla f_\rho(x_t)= \rho x_t (1-x_t)$, known as \emph{logistic map}.
    The behavior of this dynamic is highly dependent on the value of $\rho$. On the first hand, for $\rho< 3$, $\mathcal A(f_\rho, x_0)$ converges for any $x_0$. On the other hand, for almost all values of $\rho$ close enough to 4, this dynamic is chaotic.
    Note however, that for any $\rho_0<4$, there exists $\rho > \rho_0$ such that Gradient descent with step-size 1 cycles on $f_\rho$.

    Therefore, on the class $\left\{f_\rho, \rho\in [0, 4]\right\}$ we have: functions over which $\mathcal A$  converges, functions over which it diverges because it is chaotic, but also functions over which it cycles. This example thus shows that those behaviors can co-exist, but does not provide an answer to the open question.

\vspace{.2cm}
\textsc{Acknowledgments.}{\small{}The authors thank Margaux Zaffran for her feedbacks and fruitful discussions and her assistance in making plots. The work of B. Goujaud and A. Dieuleveut is partially supported by ANR-19-CHIA-0002-01/chaire SCAI, and Hi!Paris. A.~Taylor acknowledges support from the European Research Council (grant SEQUOIA 724063). This work was partly funded by the French government under management of Agence Nationale de la Recherche as part of the ``Investissements d’avenir'' program, reference ANR-19-P3IA-0001 (PRAIRIE 3IA Institute).}

\bibliographystyle{plain}
\bibliography{references}

\end{document}

%% file: header.tex
\usepackage{times}
\usepackage[utf8]{inputenc}
\usepackage[T1]{fontenc}
\usepackage{url}
\usepackage{booktabs}
\usepackage{multirow}
\usepackage{amsfonts}
\usepackage{nicefrac}
\usepackage{microtype}
\usepackage{faktor}
\usepackage{hhline}

\usepackage{thm-restate}

\usepackage{tikz}
\usepackage{mathtools, stmaryrd}
\usepackage{tocloft}
\usepackage[hidelinks]{hyperref}
\usepackage{xcolor}
\hypersetup{
    colorlinks,
    linkcolor={red!50!black},
    citecolor={blue!50!black},
    urlcolor={blue!80!black}
}
\definecolor{ao(english)}{rgb}{0.0, 0.5, 0.0}

\usepackage[algo2e, algoruled, boxed, vlined]{algorithm2e}
\usepackage{algorithm}
\usepackage{algpseudocode}

\SetKwInput{KwInput}{Input}
\SetKwInput{KwOutput}{Output}

\usepackage[noabbrev]{cleveref}

\newtheorem{Th}{Theorem}[section]

\newtheorem{Prop}[Th]{Proposition}

\newtheorem{Def}[Th]{Definition}

\newtheorem{?}[Th]{Problem}
\newtheorem{Ex}[Th]{Example}

\newenvironment{proof}[1]
{\noindent \textit{Proof.} 
}
{
$\hfill\blacksquare$}

\usepackage{amssymb}
\usepackage{pifont}

\usepackage{float}
\usepackage{wrapfig}

\usepackage{booktabs}
\usepackage{mdframed}
\newmdenv[topline=false,rightline=false]{leftbot}